\RequirePackage{fix-cm} 
\documentclass[a4paper,twoside,11pt,dvips,reqno]{amsart}

\usepackage{fixltx2e,a4wide}  
\usepackage[latin1]{inputenc}
\usepackage[T1]{fontenc}
\usepackage{amsgen}
\usepackage{amssymb}
\usepackage{amsmath}
\usepackage{mathrsfs,dsfont}
\usepackage{eucal,esint,dsfont,mathrsfs,MnSymbol,verbatim}

\usepackage{morefloats}
\usepackage{indentfirst}
\usepackage{acronym}
\usepackage{microtype}
\usepackage[labelfont={bf,sf},%
            labelsep=period,%
            justification=raggedright]{caption}
\usepackage[usenames, dvipsnames, pdf]{pstricks}
\usepackage{epsfig}
\usepackage{pst-grad} 
\usepackage{pst-plot} 

\usepackage[square,numbers]{natbib}

\newcommand{\ui}{\mathrm{i}}
\renewcommand{\Re}{\mathrm{Re}}
\renewcommand{\Im}{\mathrm{Im}}

\newcommand{\phis}{{\phi_\text{s}}}

\newcommand{\psis}{{\psi_\text{s}}}

\newcommand{\us}{{u_\text{s}}}
\newcommand{\vs}{{v_\text{s}}}

\newcommand{\zs}{{z_\text{s}}}

\newcommand{\phib}{{\phi_\text{b}}}
\newcommand{\psib}{{\psi_\text{b}}}

\newcommand{\ub}{{u_\text{b}}}
\newcommand{\vb}{{v_\text{b}}}
\newcommand{\zb}{{z_\text{b}}}

\newcommand{\nab}{\boldsymbol \nabla}
\newcommand{\scal}{\boldsymbol{\cdot}}

\newcommand{\half}{{\textstyle{1\over2}}}

\newcommand{\third}{{\textstyle{1\over3}}}

\newcommand{\sixth}{{\textstyle{1\over6}}}

\newcommand{\eqdef}{\stackrel{\text{\tiny{def}}}{=}}

\newlength{\intwidth}

\title[Explicit Dirichlet--Neumann operator]{\bf Explicit Dirichlet--Neumann operator for  
water waves}

\author[D. CLAMOND]{Didier CLAMOND}

\newcommand{\nfont}{\fontshape{n}\selectfont}
\address{({\nfont\textbf{Didier Clamond}}) Universit\'e C\^ote d'Azur, CNRS UMR 7351, 
Laboratoire J. A. Dieudonn\'e, Parc Valrose, F-06108 Nice cedex 2, France.} 
\email{didier.clamond@univ-cotedazur.fr}

\date{\today}

\begin{document}
\label{firstpage}
\maketitle

\begin{abstract}
An explicit expression for the Dirichlet--Neumann operator for surface water waves is 
presented. For non-overturning waves, but without assuming small amplitudes, the formula 
is first derived in two dimensions, subsequently extrapolated in higher dimensions and 
with a moving bottom. Although described here for water waves, this elementary approach 
could be adapted to many other problems having similar mathematical formulations.
\end{abstract}


\section{Introduction}

In this note, we consider the classical problem of gravity waves propagating at the 
(non-overturning) free surface of a homogenous non-viscous fluid in irrotational motion 
over an impermeable (uneven but non-overturning) seabed. Mathematically, in two 
dimensions without obstacles (i.e., for a simply connected fluid domain extending to 
infinity in all horizontal directions), this leads to the system of 
equations (for $x\in\mathds{R}$, $t\in\mathds{R}$ or $t\geqslant t_0$) \citep{WehausenLaitone1960}
\begin{align}
\partial_x^{\,2}\/\phi\ +\ \partial_y^{\,2}\/\phi\ 
&=\ 0 \qquad\mathrm{for}\quad -d(x)\,\leqslant\, y\,\leqslant\,\eta(x,t), \label{eqlap}\\
\partial_y\/\phi\ +\  (\partial_x\/d)\,(\partial_x\/\phi)\ 
&=\ 0 \qquad\mathrm{\,at}\qquad y\,=\,-d(x), \label{botimp}\\
\partial_y\/\phi\ -\ \partial_t\/\eta\ -\  (\partial_x\/\eta)\,(\partial_x\/\phi)\ 
&=\ 0 \qquad\mathrm{\,at}\qquad y\,=\,\eta(x,t), \label{surimp}\\
\partial_t\/\phi\ +\ g\,\eta\ +\ \half\,(\partial_x\/\phi)^{2}\ +\  
\half\,(\partial_y\/\phi)^{2}\ &=\ 0 
\qquad \mathrm{\,at}\qquad y\,=\,\eta(x,t),\label{berphi}
\end{align}
where $\phi(x,y,t)$ is a velocity potential such that $u\eqdef\partial_x\phi$ is 
the horizontal velocity and $v\eqdef\partial_y\phi$ is the vertical one, $g>0$ is 
the acceleration due to gravity (directed downward), with $(x,y)$ respectively the 
horizontal and upward-vertical Cartesian coordinates, and $t$ is the time.  
$y=\eta(x,t)$, $y=0$ and $y=-d(x)$ are, respectively, the 
equations of the free surface, of the still-water level and of the bottom; 
$h(x,t)\eqdef\eta(x,t)+d(x)$ is the total water depth. 
Physically, equation \eqref{eqlap} means that the motion is irrotational and 
isochoric, equations \eqref{botimp} and \eqref{surimp} characterise 
the impermeability of the bottom and of the free surface, while \eqref{berphi} 
expresses that the pressure at the free surface equals the constant atmospheric pressure
(set to zero without loss of generality). Capillarity and other 
surface effects can be considered but they do not affect the analysis below, 
so they are of no interest here. Also, extensions of equations \eqref{eqlap}--\eqref{berphi} 
in higher dimensions and/or moving bottoms are straightforward; these generalisations 
are considered at the end of the present paper. However, further generalisations (e.g., 
overturning surface and/or bottom, submerged obstacles, floating bodies, lateral 
solid boundaries, rough bottom) are beyond the scope of the present study; they 
require {\em ad hoc\/} investigations.  

A Dirichlet--Neumann (or Dirichlet-to-Neumann) operator (DNO) takes as input a 
function expressed at a point of the domain boundary and outputs its (outward) normal 
derivative at the same point. Here, the DNO producing the (non-unitary outgoing) normal 
derivative at the free surface is $G(\phis)\eqdef\left[\partial_y\phi-(\partial_x\eta)
(\partial_x\phi)\right]_{y=\eta}$, where $\phis(x,t)\eqdef\phi(x,\eta,t)$ denotes the velocity 
potential at the free surface. Fulfilling the Laplace equation \eqref{eqlap} and the 
bottom impermeability condition \eqref{botimp}, the DNO is a homogeneous linear function of 
$\phis$, i.e., $G(\phis)=\mathscr{G}\phis$ where $\mathscr{G}$ 
is a self-adjoint positive-definite pseudo-differential operator depending nonlinearly 
of $\eta$ and $d$ \citep{CraigSulem1993,CraigEtAl2005b}. The operator $\mathscr{G}$ is a 
fundamental mathematical object because it  `encodes' the domain geometry, the kinematic 
of the fluid motion and the bottom impermeability; moreover, it appears explicitly into 
the Hamiltonian formulation \citep{Zakharov1968} of the equations \eqref{eqlap}--\eqref{berphi}.  
Understandably, $\mathscr{G}$ has been the subject of many mathematical studies 
--- see \citet{Lannes2013} and \citet{NichollsReitich2001} for details --- and it is at 
the heart of several rigorous investigations on water waves (e.g., 
\citet{AlazardBaldi2015,AlazardEtAl2012}). 
The knowledge of the DNO mathematical features is certainly important, but its explicit 
construction is at least as important, in particular for practical applications.  
 
For flat horizontal free surface and bottom, the fluid domain is a strip and the 
DNO is easily obtained analytically, e.g., via Fourier transform. For wavy surface 
and bottom, the DNO can be constructed as a perturbation of the strip, assuming small 
amplitudes. This is the route followed in 2D by \citet{CraigSulem1993} and in 3D by 
\citet{CraigGroves1994} for flat seabeds, then extended to varying bottoms \citep{CraigEtAl2005b},  
these authors providing recurrence relations for computing the DNO to an arbitrary 
order of their  perturbative expansion.  
For small perturbations of the flat surface and seabed, other series 
representations of the DNO are available in the literature 
\citep{DommermuthYue1987,WestEtAl1987}. 
Although all these series are formally equivalent, this is not necessarily 
the case with their truncations at the same order, as outlined by 
\citet{Schaffer2008}. Moreover, such expansions are badly conditioned, so prone to 
large numerical errors and instabilities \citep{WilkeningVasan2015}. 
An explicit formulation of the DNO in expected to facilitate various reformulations for more 
efficient computations, for example, but this is not the scope of the present paper.  

The main purpose of this paper is to show how explicit Dirichlet--Neumann operators can be 
derived and, via few examples, to show their interest for analytic manipulations. Although some 
indications on potential issues and remedies with numerical computations are briefly discussed, 
it is not the purpose here to derive the most effective way to compute numerically a DNO. 

The paper is organised as follow. 
In section \ref{secDNO2D}, an explicit Dirichlet--Neumann operator is derived in two dimensions 
via rather elementary algebra. This DNO being in complex form, a real reformulation 
is introduced in section \ref{secauxrel} in order to facilitate analytical approximations. 
Some approximations for small amplitudes in finite depth and for finite amplitudes in shallow 
water are then derived in section \ref{secapp}. The DNO is extended to higher dimensions in 
section \ref{secDNO3D}, and its generalisation for moving bottoms is provided in section 
\ref{secmov}. Finally, summary and perspectives are briefly drawn in section \ref{seccon}.

\section{Two-dimensional Dirichlet--Neumann operator}\label{secDNO2D}

Let be $\psi$ the stream function harmonic conjugate of the velocity potential 
$\phi$ \citep{MilneThomson2011}. These two functions are related by the Cauchy--Riemann 
relations $\phi_x
=\psi_y=u$ and $\phi_y=-\psi_x=v$. Thus, the complex potential $f\eqdef\phi+\ui\psi$ 
is a holomorphic function of $z\eqdef x+\ui y$, with $z=\zs\eqdef x+\ui\eta$ at 
the free surface and $z=\zb\eqdef x-\ui d$ at the bottom. (As general notation,  
subscripts `s' and `b' denote quantities written, respectively, at the free surface and 
at the bottom.)
The seabed being impermeable and static, it is a streamline  
where $\psi={\psi_\text{b}}$ is constant. Without loss of generality, we then choose 
${\psi_\text{b}}=0$ for simplicity.

For any complex abscissa $z_0$, the Taylor expansion around $z_0=0$ is (omitting 
temporal dependences for brevity)
\begin{equation}\label{taylz0}
f(z-z_0)\ =\ \exp\!\left[\/-\/z_0\,\partial_z\/\right]f(z)\ \eqdef\ 
\sum_{n=0}^\infty\,
\frac{(-1)^n\,z_0^{\,n}}{n!}\,\frac{\partial^n\,f(z)}{\partial\/z^n}.
\end{equation} 
For instance, taking $z_0=\ui h=\ui(d+\eta)$, the relation \eqref{taylz0} written 
at the free surface becomes  
\begin{equation}\label{taylih}
f(\zs-\ui h)\ =\ \exp\!\left[\/-\/\ui\/h\,\partial_{\zs}\/\right]f(\zs),
\end{equation}
with the formal operator $\exp[-\ui h\partial_{\zs}]\eqdef\sum_{n=0}^\infty
(n!)^{-1}(-\ui h)^n\partial_{\zs}^{\,n}$ together with $\partial_{\zs}
\eqdef(1+\ui\eta_x)^{-1}\partial_x$, $\partial_{\zs}^{\,2}=(1+\ui\eta_x)^{-1}
\partial_x(1+\ui\eta_x)^{-1}\partial_x$, etc. (Throughout this paper, we 
use the classical convention that any operator acts on everything it multiplies on 
its right, unless parenthesis enforce otherwise.)

It should be noticed that exponents denote differential compositions, so $h^n$ is the 
$n$-th power of the function $h=d+\eta$, while $\partial_{\zs}^{\,n}$ is the $n$-th 
iteration of the differential operator $\partial_{\zs}$.
Therefore, for example, if $h$ is not constant then $h^2=h(x)^2\neq h(h(x))$, 
$h^2\,\partial_{\zs}^{\,2}
\neq(h\,\partial_{\zs})^2=h\,\partial_{\zs}\,h\,\partial_{\zs}$, 
$\exp[-\ui h\partial_{\zs}]\neq\sum_{n=0}^\infty(n!)^{-1}(-\ui h\partial_{
\zs})^n$ and the operator inverse of $\exp[-\ui h\partial_{\zs}]$ is 
{\em not\/} $\exp[\ui h\partial_{\zs}]$ (but $\exp\!\left[\ui h\partial_{\zb} 
\right]$ as shown in section \ref{secauxrel}). 

Since $\zs-\ui h=x-\ui d=\zb$ then $f(\zs-\ui h)=f(\zb)
=\phib$ is real (recall that ${\psi_\text{b}}=\Im{f_\mathrm{b}}=0$ by definition), 
while $f(\zs)=\phis+\ui\psis$ is complex. Therefore, 
the imaginary part of \eqref{taylih}, i.e., 
\begin{equation}\label{stataylIm}
0\ =\ \Re\!\left\{\/\exp\!\left[\/-\/\ui\/h\,\partial_{\zs}\/\right]\/\right\}
\psis\ +\ \Im\!\left\{\/\exp\!\left[\/-\/\ui\/h\,\partial_{\zs}\/\right]
\/\right\}\phis,
\end{equation}
yields at once
\begin{equation}\label{psisphis}
\psis\ =\ -\left(\/\Re\!\left\{\/\exp\!\left[\/-\/\ui\/h\,\partial_{\zs}
\/\right]\/\right\}\/\right)^{\!-1}\,\Im\!\left\{\/\exp\!\left[\/-\/\ui\/h\,
\partial_{\zs}\/\right]\/\right\}\phis.
\end{equation}
The equation for the free surface impermeability being $\partial_t\,\eta=\mathscr{G}\,
\phis=-\partial_x\,\psis=\vs-\us\,\partial_x\,\eta$, an explicit 
definition of the Dirichlet--Neumann operator is obtained directly from  
\eqref{psisphis} as
\begin{equation}\label{Gexp}
\mathscr{G}\ =\ \partial_x\left(\/\Re\!\left\{\/\exp\!\left[\/-\/
\ui\/h\,\partial_{\zs}\/\right]\/\right\}\/\right)^{\!-1}\,\Im\!\left\{\/
\exp\!\left[\/-\/\ui\/h\,\partial_{\zs}\/\right]\,\right\}.
\end{equation}
The formula \eqref{Gexp} provides an explicit expression for the 
DNO, i.e., $\mathscr{G}$ appears only on the left-hand side. It is the main result 
of this paper that can be generalised in higher dimensions and for moving bottoms 
(see below). It is also suitable to derive various approximations, in particular 
high-order shallow water approximations without assuming small amplitudes (see section 
\ref{secsw} below; actually, this goal was the original motivation for deriving 
\eqref{Gexp}). 

For applications, it is convenient to introduce an operator $\mathscr{J}$ such 
that $\mathscr{G}=-\partial_x\mathscr{J}\partial_x$, so
\begin{equation}\label{Jexp}
\mathscr{J}\ =\ -\left(\/\Re\!\left\{\/\exp\!\left[\/-\/
\ui\/h\,\partial_{\zs}\/\right]\/\right\}\/\right)^{\!-1}\,\Im\!\left\{\/
\exp\!\left[\/-\/\ui\/h\,\partial_{\zs}\/\right]\/\right\}\partial_x^{\,-1}.
\end{equation}
Since $\mathscr{G}$ is a self-adjoint positive-definite operator \citep{Lannes2013}, 
so is $\mathscr{J}$. Further, it is also convenient to introduce the operators 
$\mathscr{R}$ and $\mathscr{I}$ defined by
\begin{equation}
\mathscr{R}\,\eqdef\,\Re\!\left\{\/\exp\!\left[\/-\/\ui\/h\,\partial_{\zs}\/
\right]\/\right\}, \qquad 
\mathscr{I}\,\eqdef\,-\Im\!\left\{\/\exp\!\left[\/-\/\ui\/h\,\partial_{\zs}\/
\right]\/\right\}\partial_x^{\,-1}, 
\end{equation}
so $\mathscr{J}=\mathscr{R}^{-1}\mathscr{I}$. 

Although explicit, the formulae \eqref{Gexp} and \eqref{Jexp} are not quite in closed-form 
since they involve series (via the definition of the exponential operator) and operator 
inversion. Additional relations, suitable for practical applications, are then derived below.

\section{Auxiliary relations}\label{secauxrel}

With different choices of $z$ and $z_0$, the Taylor expansion \eqref{taylz0} provides 
various relations of practical interest. Several variants of \eqref{Gexp} can then be 
derived, their convenience depending on the problem at hand.

With the choice $z=\zb$ and $z_0=-\ui h=-\ui(d+\eta)$, the relation \eqref{taylz0} becomes
\begin{equation}\label{relfshfb0}
f(\zb+\ui h)\, =\, f(\zs)\, =\, \exp\!\left[\/\ui\/h\,\partial_{\zb}
\/\right]f(\zb),
\end{equation}
so a comparison with \eqref{taylih} yields at once
\begin{align}
\left(\,\exp\!\left[\/-\/\ui\/h\,\partial_{\zs}\/\right]\,\right)^{\!-1}\ 
=\ \exp\!\left[\/\ui\/h\,\partial_{\zb}\/\right] \qquad\Longleftrightarrow\qquad
\left(\/\exp\!\left[\/\ui\/h\,\partial_{\zb}\/\right]\/\right)^{\!-1}\ 
=\ \exp\!\left[\/-\/\ui\/h\,\partial_{\zs}\/\right].\label{relexpexp}
\end{align}
With this relation, the operator involving $\partial_{\zs}$ in the Dirichlet--Neumann 
operator \eqref{Gexp} can be replaced by one involving $\partial_{\zb}$. This is 
somewhat convenient in constant depth because, then, $\partial_{\zb}=\partial_x$. 
However, two operators need then to be inverted instead of one with \eqref{Gexp}, so 
further simplifications are desirable.

Taking $z=x$ together with $z_0=-\ui\eta$ and $z_0=\ui d$, \eqref{taylz0} yields
\refstepcounter{equation}
\[
f(x+\ui\eta)\, =\, f(\zs)\, =\, \exp\!\left[\/\ui\/\eta\,\partial_x\/\right]f(x),
\quad
f(x-\ui d)\, =\, f(\zb)\, =\, \exp\!\left[\/-\/\ui\/d\,\partial_x\/\right]f(x),
\eqno{(\theequation{\mathit{a},\mathit{b}})}\label{taylzsetazbd}
\]
and the elimination of $f(x)$ between these two relations, together with \eqref{taylih}, 
yields
\begin{align}
\exp\!\left[\/-\/\ui\/h\,\partial_{\zs}\/\right]&=\, 
\exp\!\left[\/-\/\ui\/d\,\partial_x\/\right]\left(\/\exp\!\left[\/\ui\/\eta\,
\partial_x\/\right]\/\right)^{-1}\, =\, \exp\!\left[\/-\/\ui\/d\,\partial_x\/\right]
\left(\/\exp\!\left[\/-\/\ui\/\eta\,\partial_x\/\right]\,\right)^\nmid\left(1+
\ui\/\eta_x\right), \label{expCR}
\end{align}
where a $\nmid$ denotes the adjoint operator.\footnote{For any complex function $\gamma$ 
of a single real variable $x$, the operator $\exp[\gamma\partial_x]\eqdef\sum_{n=0}
^\infty(n!)^{-1}\gamma^n\partial_x^{\,n}$ has for Hermitian adjoint 
$(\exp[\gamma\partial_x])^\nmid\eqdef\sum_{n=0}^\infty(n!)^{-1}\partial_x^{\,n}
(-\gamma^\ast)^n$, a star denoting the complex conjugate. We then have $(\exp[\gamma
\partial_x])^{-1}=(\exp[\gamma^\ast\partial_x])^\nmid(1+\gamma_x)$.}
We have thus relations allowing to avoid the computation of the $\partial_z$ operators, 
moreover without inversions. The operator $\mathscr{R}$ remains to be inverted, however.

For a real or complex function $\gamma$ depending on a single real variable $x$, let be 
the operators and their Hermitian adjoints
\begin{gather}
\mathscr{C}_\gamma\ \eqdef\ \sum_{n=0}^\infty\frac{(-1)^n}{(2n)!}\,\gamma^{2n}\,
\partial_x^{\,2n},\qquad
\mathscr{S}_\gamma\ \eqdef\ \sum_{n=0}^\infty\frac{(-1)^n}{(2n+1)!}\,\gamma^{2n+1}\,
\partial_x^{\,2n+1},\\
\mathscr{C}_\gamma^\nmid\ =\ \sum_{n=0}^\infty\frac{(-1)^n}{(2n)!}\,\partial_x^{\,2n}\,
{\gamma^\ast}^{2n},\qquad
\mathscr{S}_\gamma^\nmid\ =\ \sum_{n=0}^\infty\frac{(-1)^{n+1}}{(2n+1)!}\,\partial_x^{\,2n+1}\,
{\gamma^\ast}^{2n+1}.
\end{gather}
We then have 
\begin{equation}
\exp\!\left[\/-\/\ui\/d\,\partial_x\/\right]\,=\ \mathscr{C}_d\ -\ \ui\,\mathscr{S}_d,
\qquad
\left(\/\exp\!\left[\/-\/\ui\/\eta\,\partial_x\/\right]\,\right)^\nmid\ =\ 
\mathscr{C}_\eta^\nmid\ +\ \ui\,\mathscr{S}_\eta^\nmid,
\end{equation}
and the relation \eqref{expCR} is split into real and imaginary parts as
\begin{align}  
\Re\!\left\{\exp\!\left[\/-\/\ui\/h\,\partial_{\zs}\/\right]\right\}\,=\  
\mathscr{C}_d\,\mathscr{C}_\eta^\nmid\ +\ \mathscr{S}_d\,\mathscr{S}_\eta^\nmid\ -\ 
\mathscr{C}_d\,\mathscr{S}_\eta^\nmid\,\eta_x\ +\ \mathscr{S}_d\,\mathscr{C}_\eta^\nmid
\,\eta_x, \label{ReDNO2D0} \\
\Im\!\left\{\exp\!\left[\/-\/\ui\/h\,\partial_{\zs}\/\right]\right\}\,=\  
\mathscr{C}_d\,\mathscr{S}_\eta^\nmid\ -\ \mathscr{S}_d\,\mathscr{C}_\eta^\nmid\ +\ 
\mathscr{C}_d\,\mathscr{C}_\eta^\nmid\,\eta_x\ +\ \mathscr{S}_d\,\mathscr{S}_\eta^\nmid
\,\eta_x.  \label{ImDNO2D0}
\end{align}
With the operator relation $\eta_x=\partial_x\/\eta-\eta\partial_x$ (resulting from 
the Leibniz rule), we have 
\begin{align}
\mathscr{C}_\eta^\nmid\,\eta_x\ =\ \partial_x^{\,-1}\/\mathscr{S}_\eta^\nmid\,\partial_x\ 
-\ \mathscr{S}_\eta^\nmid, \qquad \mathscr{S}_\eta^\nmid\,\eta_x\ =\ \mathscr{C}_\eta^\nmid\ 
-\ \partial_x^{\,-1}\,\mathscr{C}_\eta^\nmid\,\partial_x,
\end{align}
so the relations \eqref{ReDNO2D0}--\eqref{ImDNO2D0} yield
\refstepcounter{equation}
\[
\mathscr{R}\,=\   
\mathscr{C}_d\,\partial_x^{\,-1}\,\mathscr{C}_\eta^\nmid\,\partial_x\ +\ 
\mathscr{S}_d\, \partial_x^{\,-1}\,\mathscr{S}_\eta^\nmid\,\partial_x,\qquad
\mathscr{I}\ =\ 
\mathscr{S}_d\,\partial_x^{\,-1}\,\mathscr{C}_\eta^\nmid\ -\ 
\mathscr{C}_d\,\partial_x^{\,-1}\,\mathscr{S}_\eta^\nmid. 
\eqno{(\theequation{\mathit{a},\mathit{b}})}\label{ImDNO2D}
\]
The latter relations are particularly convenient to derive analytic approximations 
and to extrapolate the DNO in higher dimensions, as shown below.

\section{Approximate Dirichlet--Neumann operators}\label{secapp}

From the explicit DNO \eqref{Gexp} and the relations derived in the previous section, 
several approximations of practical interest can be easily obtained. We consider here 
only two special cases. 

\subsection{Infinitesimal waves in arbitrary depth}\label{secCS}

Assuming that the free surface $\eta$ remains close to zero, one can formally 
expand the DNO in increasing order of nonlinearities in $\eta$ \citep{CraigSulem1993,
CraigEtAl2005b}. 
Thus, writing $\mathscr{G}=\mathscr{G}_0+\mathscr{G}_1+\mathscr{G}_2+\cdots$ 
and similarly for $\mathscr{R}$, $\mathscr{I}$ and $\mathscr{J}$, one obtains at 
once from \eqref{ImDNO2D}  
\begin{gather}
\mathscr{R}_0\, =\, \mathscr{C}_d, \qquad  
\mathscr{R}_1\, =\, -\/\mathscr{S}_d\,\eta\,\partial_x, \qquad 
\mathscr{R}_2\, =\, -\/\half\,\mathscr{C}_d\,\partial_x\,\eta^2\,\partial_x, \qquad  
\mathrm{etc.},\\
\mathscr{I}_0\, =\, \mathscr{S}_d\,\partial_x^{\,-1}, \qquad  
\mathscr{I}_1\, =\, \mathscr{C}_d\,\eta, \qquad 
\mathscr{I}_2\, =\, -\/\half\,\mathscr{S}_d\,\partial_x\,\eta^2, \qquad  
\mathrm{etc}.
\end{gather}
The relation
\begin{align}
\mathscr{R}^{-1}\ &=\,\left[\,\mathscr{R}_0\,+\,\mathscr{R}_1\,+\,\mathscr{R}_2\,
+\,\cdots\right]^{-1}\,
=\,\left[\,1\,+\,\mathscr{R}_0^{-1}\,\mathscr{R}_1\,+\,
\mathscr{R}_0^{-1}\,\mathscr{R}_2\,+\,\cdots\right]^{-1}\mathscr{R}_0^{-1}\nonumber\\
&=\,\left[\,1\,-\,\mathscr{R}_0^{-1}\,\mathscr{R}_1\,-\,
\mathscr{R}_0^{-1}\,\mathscr{R}_2\,+\,\mathscr{R}_0^{-1}\,\mathscr{R}_1\,
\mathscr{R}_0^{-1}\,\mathscr{R}_1\,+\,\cdots\right]\mathscr{R}_0^{-1},
\end{align}
then yields after some algebra
\begin{align}
\mathscr{J}_0\ &=\ \mathscr{C}_d^{-1}\,\mathscr{S}_d\,\partial_x^{\,-1}, \qquad
\mathscr{J}_1\ =\ \eta\ +\ \mathscr{J}_0\,\partial_x\,\eta\,\partial_x\/\mathscr{J}_0, \nonumber\\
\mathscr{J}_2\ &=\ \half\,\partial_x\,\eta^2\,\partial_x\/\mathscr{J}_0\ +\ 
\mathscr{J}_0\,\partial_x\,\eta\,\partial_x\,\mathscr{J}_1\ -\ 
\half\,\mathscr{J}_0\,\partial_x^{\,2}\,\eta^2, \qquad \mathrm{etc.},
\end{align}
hence
\begin{align}
\mathscr{G}_0\ &=\ -\/\partial_x\,\mathscr{C}_d^{-1}\,\mathscr{S}_d, \qquad
\mathscr{G}_1\ =\ -\/\partial_x\,\eta\,\partial_x\ -\ \mathscr{G}_0\,\eta\,
\mathscr{G}_0, \nonumber \\
\mathscr{G}_2\ &=\ \half\,\partial_x^{\,2}\,\eta^2\,\mathscr{G}_0\ 
-\ \mathscr{G}_0\,\eta\,\mathscr{G}_1
\ -\ \half\,\mathscr{G}_0\,\partial_x\,\eta^2\,\partial_x, \qquad\mathrm{etc}.
\label{GexpCS}
\end{align}
In constant depth, the expansion of \citet{CraigSulem1993} is, as expected, recovered 
introducing the operator $\mathscr{D}\eqdef\ui\partial_x$, i.e., replacing $\partial_x$ 
by $-\ui\mathscr{D}$. With a variable bottom, the expansion of \citet{CraigEtAl2005b} is 
also recovered, expect for the definition of $\mathscr{G}_0$. Indeed, \citet{CraigEtAl2005b} 
define $\mathscr{G}_0$ with an expansion for small amplitudes of the bottom corrugation 
(i.e., $\max\left|\/d(x)-\bar{d}\/\right|$ is small, where $\bar{d}$ is the mean depth), 
and they provide a recursion formula for computing this series. In \eqref{GexpCS}, 
$\mathscr{G}_0$ is defined explicitly for arbitrary (non-overturning) bottom and no 
additional expansions are required.

\subsection{Remarks}

For higher-order approximations, the recursion formula of \citet{CraigSulem1993} 
can be used verbatim with $\mathscr{G}_0$ defined here in \eqref{GexpCS}. This approach 
is convenient for the derivation of (rather low-order) analytical approximations. 
However, with numerical computations, this recursion is prone to cancelation errors 
leading to large numerical errors and instabilities \citep{WilkeningVasan2015}. This 
problem is more pronounced in higher dimensions.\footnote{W. Craig 
(2005), private communication.}  

These difficulties come mostly from the expansion of the inverse operator $\mathscr{R}^{-1}$. 
For numerical computations, this expansion should be avoided to obtain $B=\mathscr{R
}^{-1}A$ (for some functions $A$ and $B$). It is generally more efficient to solve 
$\mathscr{R}B=A$ via an iterative procedure. 
This is a similar problem as the resolution of linear systems of equations, 
for which iterative methods are often more efficient \citep{IsaacsonKeller1994}. For the DNO, 
the relation (\ref{ImDNO2D}{\it a})\footnote{See also relation \eqref{DNO3Dflat}.} shows 
that $\mathscr{R}$ behaves (roughly) like a $\cosh$-function, so $\mathscr{R}^{-1}$ 
behaves like a sech-function. The Maclaurin series of $\cosh(z)$ having an infinite 
radius of convergence, while the one of $\operatorname{sech}(z)$ convergences only 
for $|z|<\pi/2$, this provides an informal/heuristic argument showing why $B=
\mathscr{R}^{-1}A$ should not be computed but $\mathscr{R}B=A$ should be solved 
instead. With other representations (than truncated Taylor series) of $\mathscr{R}^{-1}$, 
the computation of $B=\mathscr{R}^{-1}A$ may be efficient, however.

For linear waves in the context of a highly variable bathymetry, the improvements 
of the DNO expansion proposed by \citet{AndradeNachbin2018} could be 
exploited to reformulate the explicit DNO in a more effective form for numerical 
computations. However, when speed and high numerical accuracy are required, the DNO 
perturbation expansions are not competitive (specially for steep waves) and boundary 
integral formulations should be preferred \citep{ClamondGrue2001,FructusEtAl2005,FructusGrue2007}.

\subsection{Long waves in shallow water}\label{secsw}

For long waves in shallow water, the characteristic wavelength $L_c$ is much larger 
than the characteristic depth $d_c$, so $\sigma\eqdef d_c/L_c\ll1$ is a `shallowness' 
dimensionless small parameter. The horizontal derivative $\partial_x$ is then of 
first-order in shallowness and the DNO can be expanded is power series of $\sigma$, 
without assuming small amplitude for the waves and/or for the bottom corrugation. 
Thus, we do not need to explicitly introduce scalings to asses the order of terms, 
it is sufficient to count the number of derivatives. For instance, $\partial_x^{\,3}
\eta$, $(\partial_x^{\,2}\eta)(\partial_x\eta)$ and $(\partial_x\eta)^3$ are all of 
third-order in shallowness, as well as $\partial_x^{\,3}d$, $(\partial_x^{\,2}
d)(\partial_xd)$ and $(\partial_xd)^3$. 

We then have the shallow water even terms expansions $\mathscr{R}=\mathscr{R}_0+\mathscr{R}_2
+\mathscr{R}_4+\cdots$ (and similarly for $\mathscr{I}$ and $\mathscr{J}$) so,  
from \eqref{ImDNO2D},
\begin{align}
\mathscr{R}_0\ &=\ 1,\qquad \mathscr{R}_2\ =\ -\/\half\,d^2\,\partial_x^{\,2}\ -\ 
d\,\partial_x\,\eta\,\partial_x\ -\ \half\,\partial_x\,\eta^2\,\partial_x,\qquad
\mathrm{etc.}, \\
\mathscr{I}_0\ &=\ h,\qquad  \mathscr{I}_2\ =\ -\/\sixth\,d^3\,\partial_x^{\,2}\ -\ 
\half\,d^2\,\partial_x^{\,2}\,\eta\ -\ \half\,d\,\partial_x^{\,2}\,\eta^2\ -\ 
\sixth\,\partial_x^{\,2}\,\eta^3,\qquad \mathrm{etc.}, 
\end{align}
hence, after some algebra, 
\begin{align}
\mathscr{J}_0\ &=\ h,\qquad  \mathscr{J}_2\ =\ \half\,h^2\,d_{xx}\ +\, h\,h_x\,d_x\ 
-\ h\,d_x^{\,2}\ +\ \third\,\partial_x\,h^3\,\partial_x, \qquad \mathrm{etc.}
\end{align}
Note that $\mathscr{J}_0$ and $\mathscr{J}_2$ are obviously self-adjoint, as it 
should be. 

It should be emphasised that these approximations were obtained directly from the explicit 
DNO, considering weak variations in $x$ (i.e., long waves in shallow water) but without 
assuming small amplitudes of the free surface and of the seabed (i.e., there are no restrictions 
on the magnitude of $|\eta|$ and $|d(x)-\bar{d}|$, $\bar{d}$ being the mean depth).

\section{Dirichlet--Neumann operator in higher dimensions}\label{secDNO3D}

It is rather straightforward to extrapolate the DNO given by \eqref{Gexp} to three (and 
more) spacial dimensions. In higher dimensions, the holomorphic functions cannot 
be used but series representations remain. This feature is exploited here to obtain an 
explicit expression for the DNO in an arbitrary number of dimensions.
  
With $\boldsymbol{x}=(x_1,x_2,\cdots,x_N)\in\mathds{R}^N$ referring to the `horizontal' 
coordinates, the mathematical problem is then posed in the $(N\!+\!1)$-dimensional Cartesian 
$(\boldsymbol{x},y)$-space, with $y$ the `upward-vertical' coordinate. 
Obviously, only the two-dimensional (i.e., $N=1$) and three-dimensional (i.e., $N=2$) cases 
are of physical interest for water waves. 
Let be $\boldsymbol{\nabla}\eqdef(\partial_{x_1},\cdots,\partial_{x_N})$, $\Delta\eqdef
\boldsymbol{\nabla\cdot\nabla}$ and $\mathscr{D}\eqdef\left(-\Delta\right)^{1/2}$ denote, 
respectively, the horizontal gradient, Laplacian and semi-Laplacian operators.   

The Dirichlet--Neumann operator is naturally extended in higher dimensions extrapolating the relation 
$\mathscr{G}=-\,\partial_x\mathscr{R}^{-1}\mathscr{I}\partial_x$, 
the operators $\mathscr{R}$ and $\mathscr{I}$ having to be redefined. 
In the two-dimensional case, these operators are defined via complex expressions in 
section \ref{secDNO2D}. In order to extend these operators in higher dimensions, 
one must consider their real form \eqref{ImDNO2D}, so their extrapolation is 
natural.

One-dimensional operators involving only even-order derivatives have straightforward 
extensions in higher dimensions replacing the second-order horizontal derivative 
$\partial_x^{\,2}$ by the horizontal Laplacian $\Delta$. For instance   
\begin{align}
\mathscr{C}_d\ &\mapsto\ \sum_{n=0}^\infty\frac{(-1)^n}{(2n)!}\,d^{2n}
\,\Delta^n\ \eqdef\ \cosh\!\left(\/d\,\mathscr{D}\/\right), \\
\mathscr{C}_\eta^\nmid\ &\mapsto\ \sum_{n=0}^\infty\frac{(-1)^{n}}{(2n)!}\,
\Delta^n\,\eta^{2n}\ \eqdef\ \cosh\!\left(\/\mathscr{D}\,\eta\/\right),\\
\mathscr{S}_d\,\partial_x^{\,-1}\ &\mapsto\ \sum_{n=0}^\infty\frac{(-1)^n}{(2n+1)!}
\,d^{2n+1}\,\Delta^n\ \eqdef\ \sinh\!\left(\/d\,\mathscr{D}\/\right)\mathscr{D}^{-1}, \\
\partial_x^{\,-1}\/\mathscr{S}_\eta^\nmid\ &\mapsto\ \sum_{n=0}^\infty\frac{(-1)^{n+1}}{(2n+1)!}\,
\Delta^n\,\eta^{2n+1}\ \eqdef\ -\/\mathscr{D}^{-1}\sinh\!\left(\/\mathscr{D}\,\eta\/\right).
\end{align}
It should be emphasised that, as in the one-dimensional case, the operators do not 
commute, so for example  $\cosh\!\left(\/d\/\mathscr{D}\/\right)\neq\cosh\!\left(\/
\mathscr{D}\/d\/\right)$ and $\cosh\!\left(\/d\/\mathscr{D}\/\right)^{\!-1}\neq
\operatorname{sech}\!\left(\/d\/\mathscr{D}\/\right)$, the equalities holding only 
in constant depth because then $d\/\mathscr{D}=\mathscr{D}\/d$.

A natural extension of $\mathscr{I}\partial_x$ is thus $\mathscr{I}\nab$ with
\begin{equation}\label{defI3D}
\mathscr{I}\ \mapsto\ \sinh\!\left(\/d\,\mathscr{D}\/\right)\mathscr{D}^{-1}\/
\cosh\!\left(\/\mathscr{D}\,\eta\/\right)\ +\ \cosh\!\left(\/d\,\mathscr{D}\/\right) 
\mathscr{D}^{-1}\sinh\!\left(\/\mathscr{D}\,\eta\/\right).
\end{equation}
In order to find the extension of $\partial_x\mathscr{R}^{-1}$, the operator $\mathscr{R}$ 
given by (\ref{ImDNO2D}{\it a}) is rewritten as
\begin{equation}
\mathscr{R}\ =\ \mathscr{C}_d\,\partial_x^{\,-1}\left[\,\mathscr{C}_\eta^\nmid\, 
+\,\partial_x\,\mathscr{C}_d^{\,-1}\left(\mathscr{S}_d\,\partial_x^{\,-1}\right)
\partial_x \left(\partial_x^{\,-1}\,\mathscr{S}_\eta^\nmid\right)\,\right]\partial_x.
\end{equation}
Thus, we have the natural extension
\begin{gather}
\partial_x\,\mathscr{R}^{-1}\ \mapsto\,\left[\,\cosh\!\left(\/\mathscr{D}\/\eta\/\right)\,
+\,\mathscr{G}_0\,\mathscr{D}^{-1}\sinh\!\left(\/\mathscr{D}\/\eta\/\right)\,\right]^{-1}
\nab\scal\cosh\!\left(\/d\/\mathscr{D}\/\right)^{-1}, 
\end{gather}
where $\cosh\!\left(\/d\/\mathscr{D}\/\right)^{-1}$ is the inverse operator of 
$\cosh\!\left(\/d\/\mathscr{D}\/\right)$, and
\begin{gather}
\mathscr{G}_0\ \eqdef\ -\/\nab\scal\cosh\!\left(\/d\/\mathscr{D}\/\right)^{-1}\sinh\!
\left(\/d\/\mathscr{D}\/\right)\mathscr{D}^{-1}\,\nab.\label{defG0}
\end{gather}
Therefore, the DNO becomes at once 
\begin{equation}\label{Gexp3D}
\mathscr{G}\ =\ -\left[\,\cosh\!\left(\/\mathscr{D}\/\eta\/\right)\,+\,\mathscr{G}_0
\,\mathscr{D}^{-1}\sinh\!\left(\/\mathscr{D}\/\eta\/\right)\,\right]^{-1}
\nab\scal\cosh\!\left(\/d\/\mathscr{D}\/\right)^{-1}\/\mathscr{I}\,\nab.
\end{equation}
In order to avoid misinterpretations of the formula \eqref{Gexp3D}, it is worthy to 
re-emphasise here that: (i) any operator acts on everything it multiplies on its right, 
so \eqref{Gexp3D} should be applied successively leftward 
starting from the furthest right; (ii) exponents denote operator compositions, so an exponent 
$-1$ means an operator inversion.    

Note that $\mathscr{G}\to\mathscr{G}_0$ as $\eta\to0$. Moreover, processing as in section 
\ref{secCS} for infinitesimal waves, one finds the expansion of \citet{CraigEtAl2005b}, 
except for $\mathscr{G}_0$ that is defined implicitly by \citet{CraigEtAl2005b} but 
explicitly here.

\subsection{Constant depth}

In constant depth, $d$ commuting then with both $\mathscr{D}$ and $\nab$, 
we have the simplified relation $\mathscr{G}_0=\mathscr{D}\tanh\!\left(\/d\/
\mathscr{D}\/\right)$, while $\mathscr{I}$ and $\mathscr{G}$ become (see Appendix 
\ref{appflat} for details)
\begin{align}
\mathscr{I}\ &=\ \mathscr{D}^{-1}\left[\,\sinh\!\left(\/d\/\mathscr{D}\/\right)\/
\cosh\!\left(\/\mathscr{D}\/\eta\/\right)\,+\,\cosh\!\left(\/d\/\mathscr{D}\/\right) 
\/\sinh\!\left(\/\mathscr{D}\/\eta\/\right)\right]\ =\ \mathscr{D}^{-1}\,\sinh\!
\left(\/\mathscr{D}\/h\/\right),\\
\mathscr{G}\ &=\ -\left[\,\cosh\!\left(\/d\/\mathscr{D}\/\right)\/
\cosh\!\left(\/\mathscr{D}\/\eta\/\right)\,+\,\sinh\!\left(\/d\/\mathscr{D}\/\right)
\sinh\!\left(\/\mathscr{D}\/\eta\/\right)\,\right]^{-1}
\nab\scal\mathscr{I}\,\nab\nonumber\\
&=\ -\cosh\!\left(\/\mathscr{D}\/h\/\right)^{\!-1}\,\mathscr{D}^{-1}\,\nab\scal
\sinh\!\left(\/\mathscr{D}\/h\/\right)\nab. 
\label{DNO3Dflat}
\end{align}
A better conditioned formulation, avoiding the computation of $\mathscr{D}$, is
\begin{align}
\mathscr{G}\ =\ -\left[\,\operatorname{sech}\!\left(\/d\/\mathscr{D}\/\right)
\cosh\!\left(\/\mathscr{D}\/h\/\right)\/\right]^{-1}\,\nab\scal\left[\,
\operatorname{sech}\!\left(\/d\/\mathscr{D}\/\right)\operatorname{sinhc}\!
\left(\/\mathscr{D}\/h\/\right)h\,\right]\nab, \label{DNO3Dflatbest}
\end{align} 
with 
\begin{equation}
\operatorname{sech}\!\left(\/d\/\mathscr{D}\/\right)\eqdef\,\sum_{n=0}^\infty
\frac{E_{2n}}{(2n)!}\,d^{2n}\,\mathscr{D}^{2n}, \qquad
\operatorname{sinhc}\!\left(\/\mathscr{D}\/h\/\right)\eqdef\,\sum_{n=0}^\infty
\frac{1}{(2n+1)!}\,\mathscr{D}^{2n}\,h^{2n},
\end{equation}
where $E_n$ are the Euler numbers \citep{AS65} (since $d$ and $\mathscr{D}$ commute,  
we have $\operatorname{sech}\!\left(\/d\/\mathscr{D}\/\right)=\cosh\!\left(\/d\/
\mathscr{D}\/\right)^{-1}$). The relation \eqref{DNO3Dflatbest} 
involving only even powers of $\mathscr{D}$, only Laplacian and gradient operators 
need to be evaluated, i.e., the computation of the non-local operator $\mathscr{D}$ 
can be avoided.

The DNO appearing in Hamiltonian formulations of water waves, its functional variations 
are crucial to derive the equations of motion and to investigate stability \citep{FazioliNicholls2010}. 
Thanks to the explicit DNO \eqref{DNO3Dflat}, these variations can be obtained quite 
effortlessly. Indeed, with the relations (c.f. Appendix \ref{appflat}) 
\begin{align}
\cosh\!\left(\/\mathscr{D}\/(h+\delta h)\/\right)\,&=\ 
\cosh\!\left(\/\mathscr{D}\/h\/\right)\,+\ \mathscr{D}\,\sinh\!\left(\/\mathscr{D}
\/h\/\right)\delta h\ +\ O\!\left(\delta h^2\right), \\
\sinh\!\left(\/\mathscr{D}\/(h+\delta h)\/\right)\,&=\ 
\sinh\!\left(\/\mathscr{D}\/h\/\right)\,+\ \mathscr{D}\,\cosh\!\left(\/\mathscr{D}
\/h\/\right)\delta h\ +\ O\!\left(\delta h^2\right),
\end{align}
the first variation of the DNO is obtained at once as 
\begin{align}
\mathscr{G}(h+\delta h)\ &=\ \mathscr{G}(h)\ -\ \cosh\!\left(\/\mathscr{D}\/h\/
\right)^{\!-1}\,\nab\scal\/\cosh\!\left(\/\mathscr{D}
\/h\/\right)\delta h\,\nab\nonumber\\
&\quad+\ \cosh\!\left(\/\mathscr{D}\/h\/\right)^{-1}\/\mathscr{D}\,\sinh\!\left(\/\mathscr{D}
\/h\/\right)\delta h\,\/\mathscr{G}(h)\ +\ O\!\left(\delta h^2\right),
\label{dGdh}
\end{align}
or
\begin{align}
\mathscr{G}(h+\delta h)\ &=\ \mathscr{G}(h)\ -\ \cosh\!\left(\/\mathscr{D}\/h\/
\right)^{\!-1}\,\nab\scal\/\cosh\!\left(\/\mathscr{D}
\/h\/\right)\delta h\,\nab\ +\ \mathscr{G}(h)\,\delta h\,\/\mathscr{G}(h)\nonumber\\
&\quad-\ \cosh\!\left(\/\mathscr{D}\/h\/\right)^{-1}\/\nab\scal
\cosh\!\left(\/\mathscr{D}\/h\/\right)(\nab h)\delta h\,\/\mathscr{G}(h)\ 
+\ O\!\left(\delta h^2\right).
\label{dGdhbis}
\end{align}
Similarly, higher-order functional variations of $\mathscr{G}$ can be easily obtained. 
This is one illustration of the advantage of dealing with an explicit DNO.

\subsection{Remarks}

Since the multidimensional DNO was derived extrapolating the bidimensional 
case, one can then naturally ask if \eqref{Gexp3D} is a correct expression.

First, we note that the DNO explicit expression is not unique. For instance, 
as in 2D and as suggested by the Taylor expansion around $\eta=0$, the DNO could also be 
written $\mathscr{G}=-\nab\scal\mathscr{J}\nab$ for some operator $\mathscr{J}(d,\eta, 
\mathscr{D},\nab)$ to be specified.

In 2D (i.e., for $N=1$), one can exploit the theory of holomorphic functions  
to directly check that the explicit DNO \eqref{Gexp3D} is a correct one. This procedure is 
simply the reverse of the derivations made in \S2 and \S3. This is not possible in higher 
dimension (i.e., $N>1$) because holomorphic functions cannot be used. 
The validity of \eqref{Gexp3D} was then checked expanding it ala Craig \& Sulem, 
checking that both expansions match. (This is detailed in Appendix \ref{appflat} 
for constant depth.)

\section{Moving bottom}\label{secmov}

We consider finally the generalisation of a moving bottom, i.e., $d=d(x,t)$. 
Of course, for simplicity, we begin with the two dimensional case, the 
generalisation in higher dimensions being straightforward. 

When $\partial_td\neq0$ the bottom is no longer a streamline, so the stream function is not 
zero at the seabed, i.e., ${\psi_\text{b}}={\psi_\text{b}}(x,t)\neq0$. The lower boundary 
condition \eqref{botimp} becomes $\partial_td=\partial_x\psib=-\vb-
\ub\/\partial_xd$. 
With a moving bottom, the relations \eqref{taylz0}, \eqref{taylih} and \eqref{relfshfb0} 
still hold, but \eqref{stataylIm} becomes
\begin{equation}\label{movtaylIm}
\psib\ =\ \Re\!\left\{\/\exp\!\left[\/-\/\ui\/h\,\partial_{\zs}\/\right]\/\right\}
\psis\ +\ \Im\!\left\{\/\exp\!\left[\/-\/\ui\/h\,\partial_{\zs}\/\right]
\/\right\}\phis.
\end{equation}
The condition for the bottom impermeability yielding $\psib=\partial_x^{\,-1}\/
\partial_t\/d$, the relation \eqref{movtaylIm} gives
\begin{equation}\label{psisphismov}
\psis\ =\ \Re\!\left\{\/\exp\!\left[\/-\/\ui\/h\,\partial_{\zs}\/\right]\/
\right\}^{\!-1}\left(\,\partial_x^{\,-1}\,\partial_t\,d\,-\,\Im\!\left\{\/\exp\!\left[\/-\/\ui\/
h\,\partial_{\zs}\/\right]\/\right\}\phis\,\right).
\end{equation} 
The relation \eqref{psisphismov} shows that the Dirichlet to Neumann transformation at 
the free surface is no longer a homogeneous linear function of $\phis$. The 
impermeability of the free surface is then $\partial_t\eta=G(\phis)$, the generalised 
Dirichlet--Neumann operator $G$ being
\begin{equation}
G\!\left(\phis\right)\, =\ \mathscr{G}\,\phis\ -\ \partial_x\,
\Re\!\left\{\/\exp\!\left[\/-\/\ui\/h\,\partial_{\zs}\/\right]\/\right\}^{\!-1}
\partial_x^{\,-1}\,\partial_t\,d,
\end{equation}
where $\mathscr{G}$ is given by \eqref{Gexp}. Note that $\partial_x^{\,-1}\/\partial_t\/d$ 
is not uniquely defined due to the antiderivative, unicity being enforced by the definitions 
of the mean water level and of the frame of reference. 

In higher dimension, with $\partial_x^{\,-1}=\partial_x\/\partial_x^{\,-2}\mapsto
\nab\/\Delta^{-1}=-\mathscr{D}^{-2}\/\nab$, the DNO obviously becomes
\begin{equation}\label{Gexp3Dt}
G\!\left(\phis\right)\/ =\/\left[\,\cosh\!\left(\/\mathscr{D}\,\eta\/\right)\,+\,\mathscr{G}_0
\,\mathscr{D}^{-1}\sinh\!\left(\/\mathscr{D}\,\eta\/\right)\,\right]^{-1}
\nab\scal\cosh\!\left(\/d\,\mathscr{D}\/\right)^{-1}\left(\,\mathscr{D}^{-2}\,\nab\,
\partial_t\,d\,-\/\mathscr{I}\,\nab\phis\,\right),
\end{equation}
$\mathscr{I}$ and $\mathscr{G}_0$ being defined, respectively, by \eqref{defI3D} and 
\eqref{defG0}.

\section{Discussion}\label{seccon}

Using elementary algebra, we obtained explicit formulae for the Dirichlet--Neumann 
operators involved in water wave problems. We first derived the DNO for two-dimensional 
waves over a static (uneven) bottom. We then extrapolated the formula 
to higher dimensions and generalised the formula for moving bottoms. The latter generalisation 
is interesting for its applications, such as tsunami generation \citep{Iguchi2011}, 
but also because it shows that extensions to fluids stratified in several homogeneous 
layers is possible \citep{ConstantinIvanov2019,CraigEtAl2005a}. The DNO is also used in 
some water waves problems with vorticity 
\citep{ConstantinEtAl2016b,GrovesHorn2020}, and the derivation of 
an explicit DNO for rotational waves is conceivable. 

In this note, the focus is on the DNO at the free surface assuming a given bottom shape 
and motion. Obviously, one can as easily obtain the DNO at the bottom from an assumed 
free surface that, in particular, should find applications in bottom detection from 
free surface measurements \citep{FontelosEtAl2017}.
 
The explicit DNO derived here are expressed with pseudo-differential operators 
formally defined in terms of series. Such definition supposes sufficient regularity 
of the free surface and the bottom; regularity yet to be specified by rigorous mathematical 
analysis. When these regularity conditions are not met, other more general 
representations of the operators should be used instead, such as integral formulations.
Once these operators properly defined, the explicit DNO should then be usable verbatim, 
allowing the investigation of rough bottoms and waves with angular crests, for example.  

The main purpose of this paper is to show how explicit DNO can be derived and, via examples, 
to show their interest for analytic manipulations. Although some indications on potential 
issues and remedies with numerical computations are briefly discussed, it is not the purpose 
here to derive the most effective way to compute numerically the DNO. For special functions, 
their definitions via power series are often not suitable for accurate fast computations, at 
least not in every cases and without extra knowledge (e.g., periodicity, symmetries, locations 
of singularities). The situation is similar with Dirichlet--Neumann operators defined via series, 
with the substantial extra difficulty that they involve non-commutative algebra.

Dirichlet--Neumann operators appear in many fields of research in Physics (acoustics, 
elasticity, electromagnetism, etc.) and, more generally, in the theory of partial differential 
equations. The use of a DNO is not restricted to problems involving the Laplace equation; 
it is also commonly employed in close relatives, such as the Helmholtz equation.   
The elementary formal approach presented here could then be adapted in these contexts.

\appendix

\section{Some operator relations in constant depth}\label{appflat}

In constant depth, the algebra are significantly simplified because $d$ commutes with 
both $\mathscr{D}$ and $\nab$. As mentioned at the end of section \ref{secDNO3D}, we 
then have $\mathscr{G}_0=\mathscr{D}\tanh\!\left(\/d\/\mathscr{D}\/\right)$. We also 
have, from the definition of the operators,  
\begin{align}
\cosh\!\left(\/d\/\mathscr{D}\/\right)\/\cosh\!\left(\/\mathscr{D}\/\eta\/\right)\,&=\, 
\left(\,\sum_{i=0}^\infty\,\frac{d^{2i}\,\mathscr{D}^{2i}}{(2i)!}\right)
\left(\,\sum_{j=0}^\infty\,\frac{\mathscr{D}^{2j}\,\eta^{2j}}{(2j)!}\right)\,
=\ \sum_{i,j=0}^\infty\,\frac{d^{2i}\,\mathscr{D}^{2i+2j}\,\eta^{2j}}{(2i)!\,(2j)!}\nonumber\\ 
&=\ \sum_{i=0}^\infty\,\sum_{j=0}^i\,\frac{\mathscr{D}^{2i}\,d^{2i-2j}\,\eta^{2j}}{(2j)!
\,(2i-2j)!}\ =\ \sum_{i=0}^\infty\ \sum_{j\,\text{even}}^{2i}\,\frac{\mathscr{D}^{2i}\,
d^{2i-j}\,\eta^{j}}{j!\,(2i-j)!}, \\
\sinh\!\left(\/d\/\mathscr{D}\/\right)\/\sinh\!\left(\/\mathscr{D}\/\eta\/\right)\,&=\ 
\sum_{i,j=0}^\infty\,\frac{d^{2i+1}\,\mathscr{D}^{2i+2j+2}\,\eta^{2j+1}}
{(2i+1)!\,(2j+1)!}\ =\ \sum_{i=1}^\infty\,\sum_{j=0}^{i-1}\,
\frac{\mathscr{D}^{2i}\,d^{2i-2j-1}\,\eta^{2j+1}}{(2j+1)!\,(2i-2j-1)!}\nonumber\\
&=\ \sum_{i=1}^\infty\ \sum_{j\,\text{odd}}^{2i-1}\,
\frac{\mathscr{D}^{2i}\,d^{2i-j}\,\eta^{j}}{j!\,(2i-j)!}. 
\label{taylshsh}
\end{align}
Since $|\/n!\/|=\infty$ for all negative integers $n$, the summation 
$\sum_{i=1}^\infty$ in \eqref{taylshsh} can be replaced 
by $\sum_{i=0}^\infty$. Thus, we have 
\begin{gather}
\cosh\!\left(\/d\/\mathscr{D}\/\right)\/\cosh\!\left(\/\mathscr{D}\/\eta\/\right)\,+\,
\sinh\!\left(\/d\/\mathscr{D}\/\right)\/\sinh\!\left(\/\mathscr{D}\/\eta\/\right)\,=\ 
\sum_{i=0}^\infty\frac{\mathscr{D}^{2i}}{(2i)!}\sum_{j=0}^{2i}\,\frac{(2i)!\,d^{2i-j}\,
\eta^{j}}{j!\,(2i-j)!}\nonumber\\
=\ \sum_{i=0}^\infty\,\frac{\mathscr{D}^{2i}\,(d+\eta)^{2i}}{(2i)!}\ =\ 
\cosh\!\left(\/\mathscr{D}\/h\/\right).
\end{gather}
Similarly, one can easily derive the relations
\begin{align}
\sinh\!\left(\/d\/\mathscr{D}\/\right)\/\cosh\!\left(\/\mathscr{D}\/\eta\/\right)\,+\,
\cosh\!\left(\/d\/\mathscr{D}\/\right)\/\sinh\!\left(\/\mathscr{D}\/\eta\/\right)\,&=\ 
\sinh\!\left(\/\mathscr{D}\/h\/\right),\\
\cosh\!\left(\/\eta\/\mathscr{D}\/\right)\/\sinh\!\left(\/\mathscr{D}\/d\/\right)\,+\,
\sinh\!\left(\/\eta\/\mathscr{D}\/\right)\/\cosh\!\left(\/\mathscr{D}\/d\/\right)\,&=\ 
\sinh\!\left(\/h\/\mathscr{D}\/\right),\\
\cosh\!\left(\/\eta\/\mathscr{D}\/\right)\/\cosh\!\left(\/\mathscr{D}\/d\/\right)\,+\,
\sinh\!\left(\/\eta\/\mathscr{D}\/\right)\/\sinh\!\left(\/\mathscr{D}\/d\/\right)\,&=\ 
\cosh\!\left(\/h\/\mathscr{D}\/\right),
\end{align}
and, obviously,
\begin{align}
\cosh\!\left(\/\mathscr{D}\/h\/\right)\,\pm\ \sinh\!\left(\/\mathscr{D}\/h\/\right)\,
=\ \exp\!\left(\/\pm\/\mathscr{D}\/h\/\right)
\end{align}
Note that $\mathscr{D}$ commuting with $d$, but not with $\eta$ and $h$, these relations 
are not valid for uneven bottoms and, in constant depth, $\sinh\!\left(\/\mathscr{D}\/h\/
\right)\neq\sinh\!\left(\/h\/\mathscr{D}\/\right)$ for example. However, for varying 
bottoms, similar relations can be easily obtained if $\eta$ is constant. 

Taylor expansions around $\eta=0$ yield  
\begin{align}
\sinh\!\left(\/\mathscr{D}\/h\/\right)\,&=\ \sinh\!\left(\/d\/\mathscr{D}\/\right)
\left[1+\half\/\mathscr{D}^2\/\eta^2+\cdots\right]\,
+\ \cosh\!\left(\/d\/\mathscr{D}\/\right)\left[\mathscr{D}\/\eta+\sixth\/\mathscr{D}^3
\/\eta^3+\cdots\right], \\
\cosh\!\left(\/\mathscr{D}\/h\/\right)\,&=\ \cosh\!\left(\/d\/\mathscr{D}\/\right)
\left[1+\half\/\mathscr{D}^2\/\eta^2+\cdots\right]\,
+\ \sinh\!\left(\/d\/\mathscr{D}\/\right)\left[\mathscr{D}\/\eta+\sixth\/\mathscr{D}^3
\/\eta^3+\cdots\right], 
\end{align}
hence, with $\mathscr{G}_0\eqdef\mathscr{D}\tanh\!\left(\/d\/\mathscr{D}\/\right)$, 
\begin{align}
\left(\cosh\!\left(\/\mathscr{D}\/h\/\right)\right)^{\!-1}\/ &=\,\left[\,1\,+\,\mathscr{G}_0
\,\eta\,+\,\half\,\mathscr{D}^2\,\eta^2\,+\,\sixth\,\mathscr{G}_0\,\mathscr{D}^2\,\eta^3\,+\, 
\cdots\right]^{-1}\operatorname{sech}\!\left(\/d\/\mathscr{D}\/\right)\nonumber\\
&=\,\left[\,1\,-\,\mathscr{G}_0\,\eta\,-\,\half\,\mathscr{D}^2\,\eta^2\,+\,\mathscr{G}_0\,
\eta\,\mathscr{G}_0\,\eta\,+\,\cdots\,\right]\operatorname{sech}\!\left(\/d\/\mathscr{D}\/\right),\\
\nab\scal\mathscr{D}^{-1}\sinh\!\left(\/\mathscr{D}\/h\/\right)\nab\ &=\ -
\cosh\!\left(\/d\/\mathscr{D}\/\right)\left[\,\mathscr{G}_0\,-\,\nab\scal\eta\,\nab\, 
-\,\half\,\mathscr{G}_0\,\nab\scal\eta^2\,\nab\,+\,\cdots\,\right],
\end{align}
Thus, with $\mathscr{G}$ defined in \eqref{DNO3Dflatbest}, one gets 
\begin{equation}
\mathscr{G}\ =\ \left[\,1\,+\,\mathscr{G}_0\,\eta\,+\,\half\,\mathscr{D}^2\,\eta^2\,
+\, \cdots\right]^{-1}\left[\,
\mathscr{G}_0\,-\,\nab\scal\eta\,\nab\, -\,\half\,\mathscr{G}_0\,\nab\scal\eta^2\,\nab\,
+\,\cdots\,\right].
\end{equation}
Expanding the DNO as $\mathscr{G}=\mathscr{G}_0+\mathscr{G}_1+\mathscr{G}_2+\cdots$, 
with $\mathscr{G}$ defined in \eqref{DNO3Dflat}, one obtains 
\begin{gather}
\mathscr{G}_1\, =\, -\/\mathscr{G}_0\,\eta\,\mathscr{G}_0\, -\, \nab\scal\eta\,\nab, \quad
\mathscr{G}_2\, =\,-\/\half\,\mathscr{D}^2\,\eta^2\,\mathscr{G}_0\ -\ \mathscr{G}_0\,\eta\,
\mathscr{G}_1\ -\ \half\,\mathscr{G}_0\,\nab\scal\eta^2\,\nab, \quad\text{etc.},
\end{gather}
so the expansion of \citet{CraigSulem1993} is recovered. 

Substituting $h+\delta h$ for $h$, for some small $\delta h$, we have the first-order 
Taylor expansions 
\begin{align}
\cosh\!\left(\/\mathscr{D}\/(h+\delta h)\/\right)\,&\eqdef\, 
\sum_{n=0}^{\infty}\frac{\mathscr{D}^{2n}\,(h+\delta h)^{2n}}{(2n)!}\nonumber\\
&=\ \sum_{n=0}^{\infty}\frac{\mathscr{D}^{2n}\,h^{2n}}{(2n)!}\ +\ \sum_{n=1}^{\infty}
\frac{\mathscr{D}^{2n}\,h^{2n-1}\,\delta h}{(2n-1)!}\ +\ O\!\left(\delta h^2\right)
\nonumber\\
&=\ \cosh\!\left(\/\mathscr{D}\/h\/\right)\,+\ \mathscr{D}\,\sinh\!\left(\/\mathscr{D}
\/h\/\right)\delta h\ +\ O\!\left(\delta h^2\right), \\
\sinh\!\left(\/\mathscr{D}\/(h+\delta h)\/\right)\,&=\ 
\sum_{n=0}^{\infty}\frac{\mathscr{D}^{2n+1}\,h^{2n+1}}{(2n+1)!}\ +\ \sum_{n=0}^{\infty}
\frac{\mathscr{D}^{2n+1}\,h^{2n}\,\delta h}{(2n)!}\ +\ O\!\left(\delta h^2\right)
\nonumber\\
&=\ \sinh\!\left(\/\mathscr{D}\/h\/\right)\,+\ \mathscr{D}\,\cosh\!\left(\/\mathscr{D}
\/h\/\right)\delta h\ +\ O\!\left(\delta h^2\right),
\end{align}
hence
\begin{align}
\cosh\!\left(\/\mathscr{D}\/(h+\delta h)\/\right)^{-1}\,&=\,\left[\,1\,-\,
\cosh\!\left(\/\mathscr{D}\/h\/\right)^{-1}\/\mathscr{D}\,\sinh\!\left(\/\mathscr{D}
\/h\/\right)\delta h\,\right]\cosh\!\left(\/\mathscr{D}\/h\/\right)^{-1}\ +\ 
O\!\left(\delta h^2\right).
\end{align}
We are then in position to compute explicitly the functional variations of the 
DNO.


\end{document}